\documentclass[12pt]{article}           \long\def\del#1\enddel{}
\usepackage{latexsym,cite,epic}

\let\3=\ss \catcode`\"=\active \let"=\"	\parindent=20pt	\parskip=6pt

\newfont\blackboard{msbm10 scaled 1200}   
\font\blackboards=msbm7		\font\blackboardss=msbm5
\newfam\black                   \scriptfont\black=\blackboards
\textfont\black=\blackboard     \scriptscriptfont\black=\blackboardss
\def\blackb#1{{\fam\black\relax#1}}	  
\def\IN{{\blackb N}} \def\IZ{{\blackb Z}} \def\IQ{{\blackb Q}}  
  \def\IP{{\blackb P}}

\let\l=\lambda \let\m=\mu \let\n=\nu   \let\r=\rho

  \let\G=\Gamma \let\D=\Delta
\def\0{\over }    \def\1{\vec }   \def\2{{1\over2}} \def\3{{\ss}}
\def\4{{1\over4}} \def\5{\overline }   \def\6{\partial } \def\7#1{{#1}\llap{/}}
\def\8#1{{\textstyle{#1}}}        \def\9#1{{\bf {#1}}}

\def\<{\langle } \def\>{\rangle }  
 
\def\({\left( } \def \){\right) }
 
\let\and=\wedge  

\def\mao#1{\mathop{\rm {#1}}\nolimits}

\def\beq{\begin{equation}} \def\eeq{\end{equation}} \def\eeql#1{\label{#1}\eeq}
\def\bea{\begin{eqnarray}} \def\eea{\end{eqnarray}} 
\def\beqnn#1\eeq{\[#1\]}   
\def\fnote#1#2{\begingroup\def\thefootnote{#1}\footnote{#2}
           \addtocounter{footnote}{-1}\endgroup}

\def\plb#1 #2 {Phys. Lett. {\bf B#1} #2 }
\def\phr#1 #2 {Phys. Rep. {\bf  #1} #2 } 
\def\npb#1 #2 {Nucl. Phys. {\bf B#1} #2 }
\def\aph#1 #2 {Ann. Phys. {\bf #1} #2 }  
\def\jmp#1 #2 {J. Math. Phys. {\bf #1} #2 }
\def\prd#1 #2 {Phys. Rev. {\bf D#1} #2 }
\def\prl#1 #2 {Phys. Rev. Lett. {\bf #1} #2 }
\def\rmp#1 #2 {Rev. Mod. Phys.  {\bf #1} #2 }
\def\zpc#1 #2 {Z. Phys. {\bf #1C} #2 }
\def\cmp#1 #2 {Commun. Math. Phys. {\bf #1} #2 }
\def\cqg#1 #2 {Class.Quant.Grav. {\bf #1} #2 }
\def\mpl#1 #2 {Mod. Phys. Lett. {\bf A#1} #2 }
\def\ijmp#1 #2 {Int. J. Mod. Phys. {\bf A#1} #2 }

\def\ipb{\5{\hbox{\bf 1}}} \def\ip{\hbox{\bf 1}}  
\def\q{\hbox{\bf q}}

\def\BP{\begin{picture}} \def\EP{\end{picture}}		
\def\putlin#1,#2,#3,#4,#5){\put#1,#2){\line(#3,#4){#5}}} 
\def\putvec#1,#2,#3,#4,#5){\put#1,#2){\vector(#3,#4){#5}}}
\def\putc#1)#2{\put#1){\makebox(0,0)[cc]{#2}}}
\def\BC{\begin{center}}  \def\EC{\end{center}}

\def\pmb#1{\setbox0=\hbox{${#1}$}   \kern-.025em\copy0\kern-\wd0
      \kern.05em\copy0\kern-\wd0     \kern-.025em\raise.0433em\box0 }

\def\putdot#1){\putc#1){\circle*2}}	\def\putnum#1)#2{\putc#1){\pmb{_#2}}}

\long\def\new#1\endnew{{\bf #1}}

\newfont{\XLbf}{cmbx10 scaled 2800}	\newfont{\XL}{cmr10 scaled 2600}

\begin{document}	\def\wien{TUW--96/04}			 

{\hfill  alg-geom/9603007  \vskip1pt \hfill\wien }
\vskip 19mm
\centerline{\XL Weight systems for toric Calabi--Yau varieties }
\vskip 8mm
\centerline{\XL     and reflexivity of Newton polyhedra}
\begin{center} \vskip 12mm
          Harald SKARKE\fnote{*}{e-mail: skarke@tph16.tuwien.ac.at}
\vskip 5mm
       Institut f"ur Theoretische Physik, Technische Universit"at Wien\\
       Wiedner Hauptstra\3e 8--10, A-1040 Wien, AUSTRIA

\vfill                        {\bf ABSTRACT }
\end{center}

According to a recently proposed scheme for the classification of reflexive
polyhedra, weight systems of a certain type play a prominent role.
These weight systems are classified for the cases $n=3$ and $n=4$,
corresponding to toric varieties with K3 and Calabi--Yau hypersurfaces,
respectively.
For $n=3$ we find the well known 95 weight systems corresponding to
weighted $\IP^3$'s that allow transverse polynomials, whereas for
$n=4$ there are 184026 weight systems, including the 7555 weight systems for
weighted $\IP^4$'s.
It is proven (without computer) that the Newton polyhedra corresponding
to all of these weight systems are reflexive.

\vfill \noindent \wien\\[5pt] March 1996 \vspace*{9mm}
\thispagestyle{empty} \newpage
\pagestyle{plain} 

\newpage

\section{Introduction}

A large class of Calabi--Yau manifolds can be described as hypersurfaces
defined by transverse polynomials in weighted projective 4 spaces $\IP^4$.
These spaces could be classified \cite{nms,kl94}, and the resulting list
showed a remarkable property, known as mirror symmetry: To nearly each such
variety there is a partner in the list such that the Hodge numbers
$h^{1,1}$ and $h^{2,1}$ of the two varieties are interchanged.
There are strong hints that mirror symmetry is not only a property
concerning spectra, but a full symmetry of the corresponding conformal
field theories.
Batyrev \cite{ba94} has suggested a far more powerful technique for the
construction of Calabi--Yau varieties: In his framework of reflexive
polyhedra, mirror symmetry (at the level of Hodge numbers) is manifest,
and it was checked by computer \cite{ca95,klun} that the Newton
polyhedra of all 7555 weight systems for $\IP^4$'s are reflexive, i.e.
that the older class of models is contained in this newer approach.
This makes the following goals look desirable:\\
(i) The classification of reflexive polyehedra, and\\
(ii) an explanation of the reflexivity of Newton polyhedra corresponding
to $\IP^4$'s.\\
A big step towards the solution of problem (i) was taken in ref. \cite{crp}:
There it was shown that all reflexive polyhedra are bounded by
polyhedra that can be described with the help of
certain weight systems or combinations of weight systems, and an
algorithm for the classification of these weight systems was proposed.
In the present work, I present a far more efficient algorithm and use
it to find all of the weight systems involved in the construction of
reflexive polyhedra in $n\le 4$ dimensions.
Then I show that the Newton polyhedra corresponding to any of these
weight systems (which contain the 7555 old ones as a small subset) are
reflexive, thus solving problem (ii).

There is another recent development that should be mentioned here:
It was found that, through black hole condensation, there seem to be physical
transitions between string theories compactified on different Calabi--Yau
manifolds \cite{st95,gms}, implying that there are connections between
the various moduli spaces.
In terms of toric geometry, such a transition may take place if
a reflexive polyhedron describing some Calabi--Yau hypersurface is
contained in the polyhedron describing some other Calabi--Yau hypersurface.
In this way it was shown that the moduli spaces of all Calabi--Yau
hypersurfaces of weighted $\IP^4$'s are connected \cite{cggk,acjm}.
The present work also provides a big step towards showing the connectedness
of the moduli spaces of all toric Calabi--Yau varieties: Using the fact
that the maximal Newton polyhedra corresponding to any of the weight systems
found here are reflexive and that any reflexive polyhedron is contained
in one of them, all that is left to do is to show the connectedness
of the maximal Newton polyhedra, which should be a straightforward
application of the tools developed in \cite{cggk,acjm}.

In section 2 I give a definition of reflexivity and a description of
some of the main ideas of ref. \cite{crp} used here. Then I describe
the new algorithm for the construction of weight systems required for
the classification of reflexive polyhedra and report the results
of the implementation of this algorithm on a computer.
In section 3 I give a proof that the Newton polyhedra corresponding to
any such weight system (or combination of weight systems) is
reflexive. I also give an explicit proof that the 7555 weight systems
for weighted $\IP^4$'s fall into the new set of weight systems.

\section{The classification of weight systems}

Consider a dual pair of lattices $\G\simeq \IZ^n$ and $\G^*$ and their rational
extensions $\G_\IQ$ and $\G_\IQ^*$; we denote the duality pairing
$\G_\IQ^*\times\G_\IQ\to\IQ$ by $\<\;,\;\>$.
A reflexive polytope $\D$ is an integer polytope in $\G_\IQ$
(i.e., a polytope in $\G_\IQ$ with vertices in $\G$) with exactly
one integer interior point (which we may choose to be the origin) such that
\beq
	\D^*:=\{y\in \G^*_\IQ:\<y,x\>\ge-1\:\forall x\in\D\}
\eeq
is an integer polytope in $\G^*_\IQ$.
This is equivalent to the statement that all facets of $\D$ lie on
hyperplanes of integer distance 1 to the interior point
(a lattice hyperplane $H$ has integer distance $k$ to a lattice point $P$
if there are $k-1$ lattice hyperplanes parallel to $H$ between $H$ and $P$).

Let me briefly outline the basic ideas of the algorithm proposed in \cite{crp}
for the classification of reflexive polyhedra:
Given a reflexive polytope $\D$, we look for a set of hyperplanes
$H_i$, $i=1,\cdots ,k$ carrying facets of $\D$ such that the $H_i$
define a (generically non--integer) bounded polyhedron $Q\subset\G_\IQ$.
We also assume that $Q$ is minimal in the sense that there is no
polytope with the same properties but with a smaller number of facets.
Each $H_i$ corresponds to some vertex $\5V^i$ of $\D^*$,
and $Q^*$ (the convex hull of the $\5V^i$) is an integer polytope with
the interior point of $\D^*$ in its interior.
In ref. \cite{crp} we have defined a redundant coordinate system where
the $i$'th coordinate $P^i$ of some point $P\in\G$ is given by its
integer distance to the hyperplane $H_i$ (positive on the side of $\D$).
In this way we get a natural embedding of $\G$ in $\G^k\simeq\IZ^k$.
Whenever we use this type of coordinate system, we label the interior point,
which has coordinates $(1,\cdots,1)^T$, by $\ip$.
Making use of duality and the fact that the $H_i$ have distance 1 to $\ip$,
we see that $P^i=\<\5V^i,P\>+1$.
This coordinate system has several disadvantages: We require more
coordinates than with a $\IZ^n$ description, there is an ambiguity
about the choice of $Q$, and even the lattice is not always completely
determined.
The advantages, however, seem to be greater: Our description
naturally leads to pairing matrices between vertices of $\D^*$ and $\D$
which characterise dual pairs uniquely up to the choice of some sublattice,
and even this finite ambiguity vanishes when we consider pairing
matrices for all integer points of $\D^*$ and $\D$.
In this way we avoid all the cumbersome considerations about
equivalences of polyhedra that are mapped to each other by $GL(n,\IZ)$
transformations.

Then we have shown that $Q^*$ is composed of simplices (perhaps of lower
dimension) that have $\ip$ in their interiors.
To each of these simplices there corresponds a weight system {\bf q} in the
following way:
We define the weights $q_i$ to
be the barycentric coordinates of $\ipb$ (the interior point of $\D^*$)
w.r.t. the vertices $\5V^i$ of the simplex,
i.e. $\ipb=\sum q_i\5V^i$ with $\sum q_i=1$.
The $q_i$ are positive because $\ipb$ is in the interior of the simplex
defining them.
Then $\<\ipb,P\>=0$ for $P\in\G_\IQ$ implies $\sum q_i\<\5V^i,P\>=0$,
i.e. $\sum q_i(P^i-1)=0$ and $\sum q_iP^i=1$.
The number of independent equations of this type (i.e., the number of
weight systems involved in the construction), is $k-n$.
They define an n dimensional lattice $\G^n$ with $\G\subseteq\G^n\subset\G^k$.
Now it is easy to see that $\ip$ is the only integer point in the interior
of $Q$:
For any point in the interior of $Q$ all coordinates have to be positive,
i.e. for any integer point they have to be $\ge1$.
Comparing this with $\sum q_i(P^i-1)=0$ it is clear that this can be
fulfilled only by $P^i=1$ $\forall i$.
For the construction of reflexive polyhedra we certainly need weight
systems where $\ip$ is in the interior not only of $Q$, but also in
the interior of the maximal Newton polyhedron $\D_{\rm max}=Q\cap\G^n$
defined by a weight system.
The classification algorithm proposed in \cite{crp} involved the
consideration of minimal polytopes both in $\G$ and in $\G^*$ and
the construction of pairing matrices between these polytopes.
There is, however, a far simpler way of constructing all allowed weight
systems.

The new algorithm is based on the following observation:
Assume that a weight system $q_1,\cdots,q_l$
\footnote{In this paper I always denote the dimension of $\G$ by $n$
and the number of weights in a weight system by $l$; if $\G^n$ is defined
by a single weight system, then $l=n+1$}
allows a collection of points with
coordinates $x^i$, including the interior point with $x^i=1\;\forall i$.
If these points fulfill an equation of the type $\sum_{i=1}^la_ix^i=1$
with {\bf a}$\ne${\bf q}, then the weight system must also allow at least one
point with $\sum_{i=1}^la_ix^i>1$ and at least one point with
$\sum_{i=1}^la_ix^i<1$ to ensure that {\bf 1} is really in the interior
of the maximal Newton polyhedron defined by {\bf q}.
The latter inequality is the one that we actually use for the algorithm:
Starting with the point $\ip=(1,\cdots,1)^T$, we see that unless our
weight system is $\q=(1/l,\cdots,1/l)$, there must be at least one point
with $\sum_{i=1}^lx^i<l$. For $l\le5$ there are only a few possibilities,
and after choosing some point {\bf x}${}_1$, we can look for some simple
equation fulfilled by $\ip$ and {\bf x}${}_1$ and proceed in the same way.

For $l=3$ the classification is easily carried out by hand:
Unless $\q=(1/3,1/3,1/3)$, we need at least one point with $x^1+x^2+x^3<3$.
As points where no coordinate is greater than 1 would be in conflict with
the positivity of the weight system,
we need the point $(2,0,0)^T$ (up to a permutation of indices).
Now we note that $\ip$ and $(2,0,0)^T$ both fulfill $2x^1+x^2+x^3=4$,
so $\q=(1/2,1/4,1/4)$ or we need a point with $2x^1+x^2+x^3<4$.
The only point allowed by this inequality which leads to a sensible
weight system is $(0,3,0)^T$, leading to $\q=(1/2,1/3,1/6)$.
One should note how easily we have reproduced all weight systems in
comparison with the rather lengthy analysis in ref. \cite{crp}.

For $l=4$ we can either get $\q=(1/4,1/4,1/4,1/4)$ or we need a point
with $x^1+x^2+x^3+x^4<4$. Up to permutations, all possibilities are
exhausted by $(3,0,0,0)^T$, $(2,1,0,0)^T$ and $(2,0,0,0)^T$.
For the rest of the task the computer program
requires only a few seconds. The result is a list of 99 weight systems
which still have to be checked with respect to the property that the
maximal Newton polyhedra defined by them must have $\ip$ in their interiors.

For $l=5$ similar considerations show that, unless $q_i=1/5$ for
$i=1,\cdots,5$, the weight system must allow at least one of the
points $(4,0,0,0,0)^T$, $(3,1,0,0,0)^T$, $(2,2,0,0,0)^T$, $(2,1,1,0,0)^T$,
$(3,0,0,0,0)^T$, $(2,1,0,0,0)^T$ and $(2,0,0,0,0)^T$.
Given these starting points, a C program running on an HP 735/125 required two
days of system time to produce 200653 candidates for weight systems.

The next task is to find out whether the maximal Newton polyhedra defined by
the weight systems really have $\ip$ in their interiors.
It is straightforward to construct all points allowed by some {\bf q}.
Then one could in principle construct all facets and check that
$\ip$ does not lie on one of them or on the wrong side of one of them.
I have used a different approach: Starting with $l$ points of the
maximal Newton polyhedron which are independent in $\IQ^l$, it is easy to
calculate the barycentric coordinates of $\ip$ w.r.t these points.
If all of the barycentric coordinates are positive (in this case we can
identify them with the $\5{\q}$ system introduced in \cite{crp}),
$\ip$ is in the interior of the simplex defined by the $l$ points.
If some of the barycentric coordinates
are negative, one can substitute the point corresponding to the smallest
coordinate by a point on the other side of the hyperplane defined by the
remaining points and try the same procedure again.
The same strategy can be used in cases where only $l-1$ of the starting
points are independent.
If one of the barycentric coordinates is 0 while all others are
positive, the points corresponding to the positive coordinates
define a codimension one hyperplane with $\ip$ in its interior,
so one has to check whether there is at least one point on either
side of this hyperplane.
Depending on the starting points, this strategy produced results
more or less quickly.
The best version turned out to be the one where the starting points were
determined in a fashion very similar to the algorithm that produced
the original candidates for weight systems.
Another good strategy is to use points with maximal exponents as starting
points.

It turns out that exactly 95 of the 99 weight systems for $l=4$ have the
property that $\ip$ is in the interior of the corresponding
maximal Newton polyhedron.
These are precisely the well known 95 weight systems for weighted $\IP^4$'s
that have K3 hypersurfaces \cite{reid,fl89}.

For $l=5$ the situation is completely different:
The 7555 weight systems corresponding to weighted $\IP^4$'s that allow
transverse polynomials are just a small subset of the 184026 different
weight systems whose maximal Newton polyhedra contain $\ip$.
Later I will give a proof that in arbitrary dimensions weight systems
corresponding to weighted $\IP^n$'s have the property that their
maximal Newton polyhedra contain $\ip$.
The simplest weight systems that do not correspond to weighted $\IP^4$'s
are $(1,1,1,3,4)/10$ and $(1,1,1,4,5)/12$. Note that in the latter system
the first four weights are of Fermat type, whereas the last weight
is such that no monomial of the type $X^i(X^5)^\l$, which would be necessary
for transversality, is allowed.
The corresponding maximal Newton polyhedron has vertices $V_1=(12,0,0,0,0)^T$,
$V_2=(0,12,0,0,0)^T$, $V_3=(0,0,12,0,0)^T$, $V_4=(0,0,0,3,0)^T$,
$V_5=(2,0,0,0,2)^T$, $V_6=(0,2,0,0,2)^T$ and $V_7=(0,0,2,0,2)^T$.
The facets correspond to the hyperplanes $H_i: x^i=0$ for $i=1,\cdots,5$
and $H_6: 2x^4+3x^5=6$ (spanned by the $V_j$ with $j\ge 4$).
As $\ip$ fulfils $2x^4+3x^5=5$, it has integer distance 1 to
$H_6$ and the maximal Newton polyhedron is reflexive.
Its vertex pairing matrix, i.e. the matrix ${A^i}_j=\<\5V^i,V_j\>+1$
(with $\5V^i$ corresponding to $H_i$ for $i=1,\cdots,6$) is
\beq \pmatrix{12 & 0 & 0 & 0 & 2 & 0 & 0 \cr
              0 & 12 & 0 & 0 & 0 & 2 & 0 \cr
              0 & 0 & 12 & 0 & 0 & 0 & 2 \cr
              0 & 0 & 0 & 3 & 0 & 0 & 0 \cr
              0 & 0 & 0 & 0 & 2 & 2 & 2 \cr
              6 & 6 & 6 & 0 & 0 & 0 & 0 \cr}.    \eeq
The zeroes in the matrix indicate incidence relations: ${A^i}_j=0$
means that $V_j$ lies on $H_i$.

With the help of a program that can test any weight system with respect to
the property of having $\ip$ in the interior of the maximal Newton polyhedron,
there is an easy way to check the program for the construction of all such
weight systems:
Given a positive integer $d$, one can consider all weight systems of the type
$q_i=n_i/d$ with $n_i\in\IN$ and $\sum_{i=1}^ln_i=d$ and apply the interior
point check.
I have done this up to $d=230$ for $l=5$, resulting in approximately
50000 allowed weight systems which are identical with the ones the
original program produced.

All of the weight systems found in this way might be interesting
for theoretical reasons, in particular for a better understanding of the
connection between older approaches to weighted $\IP^n$'s and the toric
framework and for dealing with the question of whether the moduli space
of all Calabi--Yau varieties allowing a description in terms of reflexive
polyhedra is connected.
For our classification program, however, we need only those weight systems
where every coordinate hyperplane is spanned by points of the maximal Newton
polyhedron.
It is easy to write a program that checks for this property.
58 of the 95 weight systems for $l=4$ pass the test (see table II in the
appendix).
For $l=5$ only approximately one fifth of the weight systems is such
that each coordinate hyperplane is spanned by points of the maximal Newton
polyhedron.
Among the 7555 weights corresponding to weighted $\IP^4$'s slightly more
than half have this property.

As an illustration for the fact that weight systems without the above
mentioned property are redundant in the classification scheme, consider
the system $(40,41,486,1134,1701)/3402$.
Its maximal Newton polyhedron $\D_{\rm max}$ is the simplex whose vertices
are the columns of the matrix
\beq \pmatrix{83 & 1 & 0 & 0 & 0 \cr
              2 & 82 & 0 & 0 & 0 \cr
              0 & 0 & 7 & 0 & 0 \cr
              0 & 0 & 0 & 3 & 0 \cr
              0 & 0 & 0 & 0 & 2 \cr}.    \eeq
The lines of this matrix correspond to the points in $\5\G$ dual to
the coordinate hyperplanes. Obviously the first two lines cannot correspond
to vertices of $\D_{\rm max}^*$. The lacking vertices of $\D_{\rm max}^*$ are
easily found to be $(84,0,0,0,0)$ and $(0,84,0,0,0)$.
Thus the vertex pairing matrix for $\D_{\rm max}^*$ and $\D_{\rm max}$ is
given by
\beq \pmatrix{84 & 0 & 0 & 0 & 0 \cr
              0 & 84 & 0 & 0 & 0 \cr
              0 & 0 & 7 & 0 & 0 \cr
              0 & 0 & 0 & 3 & 0 \cr
              0 & 0 & 0 & 0 & 2 \cr}    \eeq
which obviously corresponds to the Fermat type weight system
$(1/84,1/84,1/7,1/3,1/2)$.

We finish this section with a table of the numbers of various types
of weight systems for $l=5$.
In this table, ``span'' means the weight systems where each coordinate
hyperplane is spanned by points of the maximal Newton polyhedron, $\IP^4$ means
that the weights correspond to weighted $\IP^4$'s that allow transverse
polynomials, and in addition I have given the numbers of weight systems
containing  and not containing a weight of $1/2$.

\noindent
\begin{tabular}{||l|c|c|c|c|c|c||} \hline\hline
& $\IP^4,\;q_5=1/2$ & $\IP^4,\;q_5<1/2$ & $\IP^4$ & $q_5=1/2$ & $q_5<1/2$
                                                                  & total\\
\hline
span & 1309 & 2860 & 4169 & 14872 & 23858 & 38730 \\ \hline
total& 2390 & 5165 & 7555 & 97036 & 86990 &184026 \\
\hline\hline \end{tabular}\hfill\\[3mm]
\hfil Table I: Numbers of various types of weight systems with 5 weights

Tables III and IV in the appendix contain small sublists of the complete
list of weight systems (with small and large $d$).

\section{Some results derived without a computer}

The main aim of this section is to show that the maximal Newton polyhedra
corresponding to weights or combinations of weights constructed in the way
reported in the previous section are all reflexive.
As a prerequisite, we first need a technical lemma.

\noindent
{\bf Lemma 1:} Consider an integer pyramid $Pyr$ of height $h\ge 2$ in a
lattice
$\G\simeq\IZ^4$ and the pyramid $Pyr_{\rm double}$, which has the same peak
and the same shape as $Pyr$, but double height $2h$.
Then $Pyr_{\rm double}$ contains integer lattice points which are neither in
$Pyr$ nor in the base of $Pyr_{\rm double}$.\\[4pt]
{\it Proof:} Obviously it is sufficient to consider the case where
the base of $Pyr$ is a simplex in $\IZ^3$ (otherwise triangulate the
base and pick any simplex). Then the base of $Pyr_{\rm double}$ is a
simplex in $(2\IZ)^3$ and we may choose its vertices to be
$\1 0, 2\1e_1,2\1e_2$ and $2\1e_3$.
The peak has coordinates $(2x,2y,2z,2h)^T$.
The points in $Pyr_{\rm double}$ can be parameterized as
\beq
\l\pmatrix{2\cr 0\cr 0\cr 0\cr}+\m\pmatrix{0\cr 2\cr 0\cr 0\cr}
+\n\pmatrix{0\cr 0\cr 2\cr 0\cr}+\r\pmatrix{2x\cr 2y\cr 2z\cr 2h\cr}
\eeq
with $\l\ge0,\cdots,\r\ge 0$ and $\l+\m+\n+\r\le 1$.
With $x'=x\mao{mod}h$, $y'=y\mao{mod}h$, $z'=z\mao{mod}h$, (with
$0\le x'<h$ etc.) it is easily checked that the points
with $(\l_1,\m_1,\n_1,\r_1)=(h-x',h-y',h-z',1)/2h$
and $(\l_{h-1},\m_{h-1},\n_{h-1},\r_{h-1})=(x',y',z',h-1)/2h$
are integer points at heights 1 and $h-1$, respectively.
Clearly all parameters $\l_1, \cdots$ are positive.
With $\l_1+\m_1+\n_1+\r_1+\l_{h-1}+\m_{h-1}+\n_{h-1}+\r_{h-1}=2$,
at least one of the inequalities $\l_1+\m_1+\n_1+\r_1\le 1$
and $\l_{h-1}+\m_{h-1}+\n_{h-1}+\r_{h-1}\le 1$ must be fulfilled.
This means that at least one of the two points is a point of
$Pyr_{\rm double}$, and because of the height it is neither in $Pyr$ nor
in the base.
\hfill$\Box$\\[4pt]
{\bf Remarks:} The same proof works for lattices $\IZ^n$
with $n<4$.
For $n=5$ there is the following counterexample:
Let the base again be given by $\1 0$ and $2\1e_i$ and the peak by
$(2,2,2,2,4)^T$. With $h=2$, a point fulfilling the criterion of the
lemma would have to be at height 1. With the same ansatz as in the proof,
we would have $\r=1/4$ and $\l\ge 0, \m\ge 0, \cdots$ would have to be
at least $1/4$, resulting in $\l+\cdots+\r\ge5/4$ and a point outside
$Pyr_{\rm double}$.
Thus there is no integer point between the bases of $Pyr$ and
$Pyr_{\rm double}$.

\noindent
{\bf Theorem:} Four or lower dimensional maximal Newton polyhedra with $\ip$
in their interior are reflexive.\\[4pt]
{\it Proof:} Consider a collection of points in a maximal Newton polyhedron
$\D_{\rm max}$ spanning a hyperplane at distance $h\ge 2$ from $\ip$.
We take these points to define the base of the pyramid $Pyr$ of lemma 1 and
$\ip$ as the peak.
Then $Pyr_{\rm double}$ lies in
$(\D_{\rm max})_{\rm double}\subseteq \{x\in\G^n:x^i\ge -1\}$.
Only the base of $Pyr_{\rm double}$ can intersect with the boundary of
$(\D_{\rm max})_{\rm double}$, so the integer points of $Pyr_{\rm double}$
that are not in the base must have nonnegative coordinates, i.e. they must
be in $\D_{\rm max}$.
Thus the lemma ensures that there are points in the cone defined by $Pyr$
which are outside of $Pyr$, but within $\D_{\rm max}$. This means that the
hyperplane
defined by the base of $Pyr$ is not a bounding hyperplane of $\D_{\rm max}$.
Therefore every bounding hyperplane of $\D_{\rm max}$ must be at distance 1,
i.e. $\D_{\rm max}$ is reflexive.\hfill$\Box$\\[4pt]
{\bf Remark:} The theorem holds not only for maximal Newton polyhedra defined
by a single weight system with $l=n+1$, as mainly considered in this paper,
but also for maximal Newton polyhedra defined by several weight systems with
$l<n+1$ involved in the classification scheme of \cite{crp}.\\[4pt]
{\bf Examples:} The weight
system $q_i=1/5$, $i=1,\cdots,5$ contains the points $(2,0,0,0,3)^T$,
$(0,2,0,0,3)^T$, $(0,0,2,0,3)^T$, $(0,0,0,2,3)^T$ defining the
hyperplane $x^5=3$ (at distance 2 to $\ip$). The base of $Pyr_{\rm double}$
is the convex hull of the points $(3,-1,-1,-1,5)^T$, $(-1,3,-1,-1,5)^T$,
$(-1,-1,3,-1,5)^T$, $(-1,-1,-1,3,5)^T$, and the ``height one'' points
$(1,0,0,0,4)^T$, $(0,1,0,0,4)^T$, $(0,0,1,0,4)^T$, $(0,0,0,1,4)^T$
ensure that $x^5=3$ is not a bounding hyperplane. In the same way the
hyperplane $x^5=4$ (at distance 3 to $\ip$), with the points
$(1,0,0,0,4)^T$, $(0,1,0,0,4)^T$, $(0,0,1,0,4)^T$, $(0,0,0,1,4)^T$,
gives rise to a pyramid of height 6 with base points
$(1,-1,-1,-1,7)^T$, $(-1,1,-1,-1,7)^T$, $(-1,-1,1,-1,7)^T$, $(-1,-1,-1,1,7)^T$.
This time the integer point we are looking for is $(0,0,0,0,5)^T$.

The fact that maximal Newton polyhedra of weighted $\IP^4$'s are reflexive
has been
known for some time due to explicit computer calculations \cite{ca95,klun}.
In order to rederive this result without the help of a computer, we still
have to show that the list of 7555 weights for $\IP^4$'s is contained
in our complete list of weights whose maximal Newton polyhedra have $\ip$
in the interior.
The analogous statement holds in any dimension and also for
abelian orbifolds (with the transversality condition applied to
polynomials that are invariant under the twist group); although it looks
quite obvious to anyone who has worked with weighted $\IP^n$'s for
some time, the proof turns out to be rather technical.

\noindent
{\bf Lemma 2:} Maximal Newton polyhedra corresponding to weighted $\IP^n$'s
or abelian orbifolds of weighted $\IP^n$'s that allow transverse polynomials
have $\ip$ in their interiors.\\[4pt]
{\it Proof:} A transverse polynomial contains monomials of the type
$M^i=(X^i)^{a_i}$ or $M^i=(X^i)^{a_i}X^j$ (with $a_i\ge 2$) for each $i$.
These monomials define points which can be arranged in the matrix
\beq A=\pmatrix{a_1&x&x&x&\cdots\cr
                x&a_2&x&x&\cdots\cr
                x&x&a_3&x&\cdots\cr
                \cdots&&&&\cr},     \eeq
where in each column at most one of the $x$'s can be 1, whereas all others
are zero.
Let us first see that $A$ is regular: Assuming it were singular,
we could find a nontrivial vanishing linear combination of its lines, i.e.
we would have $\1\l\ne0$ with $\sum_i \l_iA^i{}_j=0$.
The specific form of our matrix implies that $\l_i=0$ if $M^i=(X^i)^{a_i}$
and $a_i\l_i+\l_{p(i)}=0$ if $M^i=(X^i)^{a_i}X^{p(i)}$.
Iterating this, we get $a_ia_{p(i)}\l_i-\l_{p(p(i))}=0$ etc.
At some point either $\l_{p(\cdots(i)\cdots)}=0$ or $p(\cdots(i)\cdots)=i$,
showing that indeed $\l_i=0$.
Thus $A$ is regular and we can solve $\sum_i A^i{}_j\5q^j=1$ for $\5q^j$.
Then $\sum_j\5q^j=\sum_{i,j}(A^{-1})^j{}_i=\sum_iq_i=1$, showing that the
$\5q^j$ are the barycentric coordinates of $\ip$ with respect to the
columns of $A$. 
If all of the $\5q^j$ are positive, $\ip$ is in the interior of the simplex
defined by the columns of $A$.
Let us now assume that not all of the $\5q^j$ are positive:
Let $\5q^j\le 0$ for $j\in I_-$ and $\5q^j>0$ for $j\in I_+$, with
$I_-\cup I_+=\{1,\cdots,n\}$.
Now sum the equations defining the $\5q^j$ over $i\in I_-$ to get
\beq   \sum_{i\in I_-}\sum_j A^i{}_j\5q^j=|I_-|   \eeql{grausgl}
and split $\sum_j$ in $\sum_{j\in I_-}+\sum_{j\in I_+}$.
Then
\beq   \sum_{i\in I_-}\sum_{j\in I_-} A^i{}_j\5q^j\le2\sum_{j\in I_-}\5q^j\eeq
because the diagonal elements are involved and
\beq   \sum_{i\in I_-}\sum_{j\in I_+} A^i{}_j\5q^j\le\sum_{j\in I_+}\5q^j \eeq
because each column contains at most a single 1.
Thus the l.h.s. of eq. (\ref{grausgl}) fulfils
\beq   \sum_{i\in I_-}\sum_j A^i{}_j\5q^j\le
     2\sum_{j\in I_-}\5q^j+\sum_{j\in I_+}\5q^j=1+\sum_{j\in I_-}\5q^j\le 1
\eeq
with equality iff
$\sum_{j\in I_-}\5q^j=0$ and $\sum_{i\in I_-}\sum_{j\in I_+} A^i{}_j\5q^j=1$.
Therefore $I_-$ can have at most 1 element corresponding to some $\5q^j=0$, and
the corresponding line may contain no 0, i.e. up to permutations our matrix is
\beq A=\pmatrix{a_1&1&1&1&\cdots\cr
                0&a_2&0&0&\cdots\cr
                0&0&a_3&0&\cdots\cr
                \cdots&&&&\cr}.     \eeq
$\ip$ is in the interior of the $n-1$ dimensional simplex defined by
all columns except the first.
This simplex lies on the hyperplane $x^1=1$.
The first point is above this hyperplane, and the transversality condition
ensures that there is also a point with $x^1=0$, i.e. a point below
this hyperplane.
Therefore $\ip$ is again in the interior of the maximal Newton polyhedron.
\hfill$\Box$\\

\bigskip

{\it Acknowledgements.}  I would like to thank Max Kreuzer for a long
collaboration to which I owe my interest in reflexive polyhedra, and also
for helpful discussions in the context of the present work.
This work is supported
by the {\it Austrian National Bank} under grant number 5674.

\vfill\eject
\noindent
{\bf\Large Appendix: Various Tables}

In table II I list the 95 weight systems for $l=4$.
The last column indicates whether a weight system has the property that
coordinate hyperplanes are spanned by points of the maximal Newton polyhedron.

Table III contains all weight systems for $l=5$ with $d\le20$, whereas
table IV contains the weight systems with the largest values of $d$ for
the cases $q_5=1/2$ and $q_5<1/2$.
The last columns indicate whether a weight system corresponds to a weighted
$\IP^4$ that allows transverse polynomials.
All weights in table III have the property that
coordinate hyperplanes are spanned by points of the maximal Newton polyhedron,
whereas none of the weight systems in table IV fulfil this criterion.

\bigskip

\noindent
\begin{tabular}{||rrrr|r|r|} \hline\hline
$n_1$  & $n_2$  & $n_3$  & $n_4$  & d & span  \\ \hline
1 & 1 & 1 & 1 & 4 & y\\
1 & 1 & 1 & 2 & 5 & y\\
1 & 1 & 2 & 2 & 6 & y\\
1 & 1 & 1 & 3 & 6 & y\\
1 & 1 & 2 & 3 & 7 & y\\
1 & 2 & 2 & 3 & 8 & y\\
1 & 1 & 2 & 4 & 8 & y\\
1 & 2 & 3 & 3 & 9 & y\\
1 & 1 & 3 & 4 & 9 & y\\
1 & 2 & 3 & 4 & 10 & y\\
1 & 2 & 2 & 5 & 10 & y\\
1 & 1 & 3 & 5 & 10 & y\\
1 & 2 & 3 & 5 & 11 & y\\
2 & 3 & 3 & 4 & 12 & y\\
1 & 3 & 4 & 4 & 12 & y\\
2 & 2 & 3 & 5 & 12 & y\\
1 & 2 & 4 & 5 & 12 & y\\
1 & 2 & 3 & 6 & 12 & y\\
1 & 1 & 4 & 6 & 12 & y\\
1 & 3 & 4 & 5 & 13 & y\\
2 & 3 & 4 & 5 & 14 & y\\
2 & 2 & 3 & 7 & 14 & y\\
1 & 2 & 4 & 7 & 14 & y\\
3 & 3 & 4 & 5 & 15 & y\\
2 & 3 & 5 & 5 & 15 & y\\
1 & 3 & 5 & 6 & 15 & y\\
1 & 3 & 4 & 7 & 15 & y\\
1 & 2 & 5 & 7 & 15 & y\\
1 & 4 & 5 & 6 & 16 & y\\
2 & 3 & 4 & 7 & 16 & y\\
1 & 3 & 4 & 8 & 16 & y\\
1 & 2 & 5 & 8 & 16 & y\\
\hline\hline \end{tabular}\hfil
\hbox{\begin{tabular}{|rrrr|r|r|} \hline\hline
$n_1$  & $n_2$  & $n_3$  & $n_4$  & d & span  \\ \hline
2 & 3 & 5 & 7 & 17 & y\\
3 & 4 & 5 & 6 & 18 & y\\
1 & 4 & 6 & 7 & 18 & y\\
2 & 3 & 5 & 8 & 18 & y\\
2 & 3 & 4 & 9 & 18 & y\\
1 & 3 & 5 & 9 & 18 & y\\
1 & 2 & 6 & 9 & 18 & y\\
3 & 4 & 5 & 7 & 19 & y\\
2 & 5 & 6 & 7 & 20 & n\\
3 & 4 & 5 & 8 & 20 & y\\
2 & 4 & 5 & 9 & 20 & n\\
2 & 3 & 5 & 10 & 20 & y\\
1 & 4 & 5 & 10 & 20 & y\\
3 & 5 & 6 & 7 & 21 & n\\
1 & 5 & 7 & 8 & 21 & n\\
2 & 3 & 7 & 9 & 21 & y\\
1 & 3 & 7 & 10 & 21 & n\\
2 & 4 & 5 & 11 & 22 & n\\
1 & 4 & 6 & 11 & 22 & y\\
1 & 3 & 7 & 11 & 22 & n\\
3 & 6 & 7 & 8 & 24 & n\\
4 & 5 & 6 & 9 & 24 & y\\
1 & 6 & 8 & 9 & 24 & y\\
3 & 4 & 7 & 10 & 24 & y\\
2 & 3 & 8 & 11 & 24 & n\\
3 & 4 & 5 & 12 & 24 & y\\
2 & 3 & 7 & 12 & 24 & y\\
1 & 3 & 8 & 12 & 24 & y\\
4 & 5 & 7 & 9 & 25 & n\\
2 & 5 & 6 & 13 & 26 & n\\
1 & 5 & 7 & 13 & 26 & n\\
2 & 3 & 8 & 13 & 26 & n\\
\hline\hline \end{tabular}}\hfil
\hbox{\begin{tabular}{|rrrr|r|r||} \hline\hline
$n_1$  & $n_2$  & $n_3$  & $n_4$  & d & span  \\ \hline
5 & 6 & 7 & 9 & 27 & n\\
2 & 5 & 9 & 11 & 27 & n\\
4 & 6 & 7 & 11 & 28 & n\\
3 & 4 & 7 & 14 & 28 & y\\
1 & 4 & 9 & 14 & 28 & n\\
5 & 6 & 8 & 11 & 30 & n\\
3 & 4 & 10 & 13 & 30 & n\\
4 & 5 & 6 & 15 & 30 & y\\
2 & 6 & 7 & 15 & 30 & n\\
1 & 6 & 8 & 15 & 30 & y\\
2 & 3 & 10 & 15 & 30 & y\\
1 & 4 & 10 & 15 & 30 & y\\
4 & 5 & 7 & 16 & 32 & n\\
2 & 5 & 9 & 16 & 32 & n\\
3 & 5 & 11 & 14 & 33 & n\\
4 & 6 & 7 & 17 & 34 & n\\
3 & 4 & 10 & 17 & 34 & n\\
7 & 8 & 9 & 12 & 36 & n\\
3 & 4 & 11 & 18 & 36 & n\\
1 & 5 & 12 & 18 & 36 & n\\
5 & 6 & 8 & 19 & 38 & n\\
3 & 5 & 11 & 19 & 38 & n\\
5 & 7 & 8 & 20 & 40 & n\\
3 & 4 & 14 & 21 & 42 & y\\
2 & 5 & 14 & 21 & 42 & n\\
1 & 6 & 14 & 21 & 42 & y\\
4 & 5 & 13 & 22 & 44 & n\\
3 & 5 & 16 & 24 & 48 & n\\
7 & 8 & 10 & 25 & 50 & n\\
4 & 5 & 18 & 27 & 54 & n\\
5 & 6 & 22 & 33 & 66 & n\\
 &  &  &  &  & \\
\hline\hline \end{tabular}}\hfill\\[3mm]
\hfil Table II: Weights for $l=4$

\noindent
\begin{tabular}{||rrrrr|r|r||} \hline\hline
$n_1$  & $n_2$  & $n_3$  & $n_4$  & $n_5$  & d & $\IP^4$  \\ \hline
1 & 1 & 1 & 1 & 1 & 5 & y\\
1 & 1 & 1 & 1 & 2 & 6 & y\\
1 & 1 & 1 & 2 & 2 & 7 & y\\
1 & 1 & 1 & 1 & 3 & 7 & y\\
1 & 1 & 2 & 2 & 2 & 8 & y\\
1 & 1 & 1 & 2 & 3 & 8 & y\\
1 & 1 & 1 & 1 & 4 & 8 & y\\
1 & 1 & 2 & 2 & 3 & 9 & y\\
1 & 1 & 1 & 3 & 3 & 9 & y\\
1 & 1 & 1 & 2 & 4 & 9 & y\\
1 & 2 & 2 & 2 & 3 & 10 & y\\
1 & 1 & 2 & 3 & 3 & 10 & y\\
1 & 1 & 2 & 2 & 4 & 10 & y\\
1 & 1 & 1 & 3 & 4 & 10 & n\\
1 & 1 & 1 & 2 & 5 & 10 & y\\
1 & 2 & 2 & 3 & 3 & 11 & y\\
1 & 1 & 2 & 3 & 4 & 11 & y\\
1 & 1 & 2 & 2 & 5 & 11 & y\\
1 & 1 & 1 & 3 & 5 & 11 & y\\
2 & 2 & 2 & 3 & 3 & 12 & y\\
1 & 2 & 3 & 3 & 3 & 12 & y\\
1 & 2 & 2 & 3 & 4 & 12 & y\\
1 & 1 & 3 & 3 & 4 & 12 & y\\
1 & 1 & 2 & 4 & 4 & 12 & y\\
1 & 2 & 2 & 2 & 5 & 12 & y\\
1 & 1 & 2 & 3 & 5 & 12 & y\\
1 & 1 & 1 & 4 & 5 & 12 & n\\
1 & 1 & 2 & 2 & 6 & 12 & y\\
1 & 1 & 1 & 3 & 6 & 12 & y\\
1 & 2 & 3 & 3 & 4 & 13 & y\\
1 & 1 & 3 & 4 & 4 & 13 & y\\
1 & 2 & 2 & 3 & 5 & 13 & y\\
1 & 1 & 3 & 3 & 5 & 13 & y\\
1 & 1 & 2 & 4 & 5 & 13 & n\\
1 & 1 & 2 & 3 & 6 & 13 & y\\
1 & 1 & 1 & 4 & 6 & 13 & y\\
2 & 2 & 3 & 3 & 4 & 14 & y\\
1 & 2 & 3 & 4 & 4 & 14 & n\\
2 & 2 & 2 & 3 & 5 & 14 & n\\
1 & 2 & 3 & 3 & 5 & 14 & n\\
1 & 2 & 2 & 4 & 5 & 14 & y\\
1 & 1 & 3 & 4 & 5 & 14 & n\\
1 & 2 & 2 & 3 & 6 & 14 & y\\
1 & 1 & 2 & 4 & 6 & 14 & y\\
1 & 2 & 2 & 2 & 7 & 14 & y\\
\hline\hline \end{tabular}\hfil
\hbox{\begin{tabular}{||rrrrr|r|r||} \hline\hline
$n_1$  & $n_2$  & $n_3$  & $n_4$  & $n_5$  & d & $\IP^4$  \\ \hline
1 & 1 & 2 & 3 & 7 & 14 & y\\
1 & 1 & 1 & 4 & 7 & 14 & n\\
2 & 3 & 3 & 3 & 4 & 15 & y\\
1 & 3 & 3 & 4 & 4 & 15 & y\\
2 & 2 & 3 & 3 & 5 & 15 & y\\
1 & 3 & 3 & 3 & 5 & 15 & y\\
1 & 2 & 3 & 4 & 5 & 15 & y\\
1 & 2 & 2 & 5 & 5 & 15 & y\\
1 & 1 & 3 & 5 & 5 & 15 & y\\
1 & 2 & 3 & 3 & 6 & 15 & y\\
1 & 1 & 3 & 4 & 6 & 15 & y\\
1 & 1 & 2 & 5 & 6 & 15 & n\\
1 & 2 & 2 & 3 & 7 & 15 & y\\
1 & 1 & 3 & 3 & 7 & 15 & y\\
1 & 1 & 2 & 4 & 7 & 15 & y\\
1 & 1 & 1 & 5 & 7 & 15 & y\\
2 & 3 & 3 & 4 & 4 & 16 & y\\
1 & 3 & 4 & 4 & 4 & 16 & y\\
2 & 2 & 3 & 4 & 5 & 16 & n\\
1 & 3 & 3 & 4 & 5 & 16 & y\\
1 & 2 & 4 & 4 & 5 & 16 & y\\
1 & 2 & 3 & 5 & 5 & 16 & n\\
1 & 1 & 4 & 5 & 5 & 16 & y\\
1 & 2 & 3 & 4 & 6 & 16 & y\\
1 & 1 & 4 & 4 & 6 & 16 & y\\
1 & 2 & 2 & 5 & 6 & 16 & n\\
1 & 1 & 3 & 5 & 6 & 16 & n\\
2 & 2 & 2 & 3 & 7 & 16 & y\\
1 & 2 & 3 & 3 & 7 & 16 & y\\
1 & 2 & 2 & 4 & 7 & 16 & y\\
1 & 1 & 3 & 4 & 7 & 16 & n\\
1 & 1 & 2 & 5 & 7 & 16 & y\\
1 & 2 & 2 & 3 & 8 & 16 & y\\
1 & 1 & 3 & 3 & 8 & 16 & y\\
1 & 1 & 2 & 4 & 8 & 16 & y\\
1 & 1 & 1 & 5 & 8 & 16 & y\\
2 & 3 & 3 & 4 & 5 & 17 & y\\
1 & 3 & 4 & 4 & 5 & 17 & n\\
2 & 2 & 3 & 5 & 5 & 17 & y\\
1 & 2 & 4 & 5 & 5 & 17 & n\\
1 & 2 & 3 & 5 & 6 & 17 & y\\
1 & 1 & 4 & 5 & 6 & 17 & n\\
2 & 2 & 3 & 3 & 7 & 17 & y\\
1 & 2 & 3 & 4 & 7 & 17 & y\\
1 & 2 & 2 & 5 & 7 & 17 & n\\
\hline\hline \end{tabular}}\hfil\\
\hbox{\begin{tabular}{||rrrrr|r|r||} \hline\hline
$n_1$  & $n_2$  & $n_3$  & $n_4$  & $n_5$  & d & $\IP^4$  \\ \hline
1 & 1 & 3 & 5 & 7 & 17 & y\\
1 & 2 & 3 & 3 & 8 & 17 & y\\
1 & 1 & 3 & 4 & 8 & 17 & y\\
1 & 1 & 2 & 5 & 8 & 17 & y\\
3 & 3 & 3 & 4 & 5 & 18 & n\\
2 & 3 & 4 & 4 & 5 & 18 & n\\
2 & 3 & 3 & 5 & 5 & 18 & y\\
1 & 3 & 4 & 5 & 5 & 18 & n\\
2 & 3 & 3 & 4 & 6 & 18 & y\\
1 & 3 & 4 & 4 & 6 & 18 & n\\
2 & 2 & 3 & 5 & 6 & 18 & y\\
1 & 3 & 3 & 5 & 6 & 18 & y\\
1 & 2 & 4 & 5 & 6 & 18 & n\\
1 & 2 & 3 & 6 & 6 & 18 & y\\
1 & 1 & 4 & 6 & 6 & 18 & y\\
2 & 2 & 3 & 4 & 7 & 18 & y\\
1 & 3 & 3 & 4 & 7 & 18 & n\\
1 & 2 & 4 & 4 & 7 & 18 & n\\
1 & 2 & 3 & 5 & 7 & 18 & n\\
1 & 1 & 4 & 5 & 7 & 18 & n\\
1 & 2 & 2 & 6 & 7 & 18 & n\\
1 & 1 & 3 & 6 & 7 & 18 & n\\
2 & 2 & 3 & 3 & 8 & 18 & y\\
1 & 2 & 3 & 4 & 8 & 18 & n\\
1 & 2 & 2 & 5 & 8 & 18 & y\\
1 & 1 & 3 & 5 & 8 & 18 & n\\
1 & 1 & 2 & 6 & 8 & 18 & y\\
2 & 2 & 2 & 3 & 9 & 18 & y\\
1 & 2 & 3 & 3 & 9 & 18 & y\\
1 & 2 & 2 & 4 & 9 & 18 & y\\
1 & 1 & 3 & 4 & 9 & 18 & n\\
1 & 1 & 2 & 5 & 9 & 18 & n\\
1 & 1 & 1 & 6 & 9 & 18 & y\\
3 & 3 & 4 & 4 & 5 & 19 & y\\
2 & 3 & 4 & 5 & 5 & 19 & n\\
1 & 3 & 4 & 5 & 6 & 19 & y\\
2 & 3 & 3 & 4 & 7 & 19 & n\\
1 & 3 & 4 & 4 & 7 & 19 & n\\
2 & 2 & 3 & 5 & 7 & 19 & n\\
1 & 3 & 3 & 5 & 7 & 19 & n\\
1 & 2 & 4 & 5 & 7 & 19 & y\\
1 & 2 & 3 & 6 & 7 & 19 & n\\
1 & 1 & 4 & 6 & 7 & 19 & n\\
1 & 3 & 3 & 4 & 8 & 19 & y\\
\hline\hline \end{tabular}}\hfil
\hbox{\begin{tabular}{||rrrrr|r|r||} \hline\hline
$n_1$  & $n_2$  & $n_3$  & $n_4$  & $n_5$  & d & $\IP^4$  \\ \hline
1 & 2 & 3 & 5 & 8 & 19 & n\\
1 & 1 & 3 & 6 & 8 & 19 & y\\
1 & 2 & 3 & 4 & 9 & 19 & y\\
1 & 2 & 2 & 5 & 9 & 19 & y\\
1 & 1 & 3 & 5 & 9 & 19 & y\\
1 & 1 & 2 & 6 & 9 & 19 & y\\
3 & 4 & 4 & 4 & 5 & 20 & y\\
3 & 3 & 4 & 5 & 5 & 20 & y\\
2 & 4 & 4 & 5 & 5 & 20 & y\\
2 & 3 & 5 & 5 & 5 & 20 & y\\
1 & 4 & 5 & 5 & 5 & 20 & y\\
2 & 3 & 4 & 5 & 6 & 20 & y\\
1 & 4 & 4 & 5 & 6 & 20 & n\\
2 & 2 & 5 & 5 & 6 & 20 & y\\
1 & 3 & 5 & 5 & 6 & 20 & n\\
1 & 2 & 5 & 6 & 6 & 20 & n\\
2 & 3 & 4 & 4 & 7 & 20 & n\\
2 & 3 & 3 & 5 & 7 & 20 & n\\
2 & 2 & 4 & 5 & 7 & 20 & n\\
1 & 3 & 4 & 5 & 7 & 20 & n\\
1 & 2 & 5 & 5 & 7 & 20 & n\\
2 & 2 & 3 & 6 & 7 & 20 & y\\
1 & 2 & 4 & 6 & 7 & 20 & y\\
1 & 1 & 5 & 6 & 7 & 20 & n\\
2 & 3 & 3 & 4 & 8 & 20 & y\\
1 & 3 & 4 & 4 & 8 & 20 & y\\
2 & 2 & 3 & 5 & 8 & 20 & n\\
1 & 3 & 3 & 5 & 8 & 20 & n\\
1 & 2 & 4 & 5 & 8 & 20 & y\\
1 & 2 & 3 & 6 & 8 & 20 & n\\
1 & 1 & 4 & 6 & 8 & 20 & y\\
2 & 2 & 3 & 4 & 9 & 20 & y\\
1 & 2 & 4 & 4 & 9 & 20 & y\\
2 & 2 & 2 & 5 & 9 & 20 & y\\
1 & 2 & 3 & 5 & 9 & 20 & y\\
1 & 1 & 4 & 5 & 9 & 20 & n\\
1 & 2 & 2 & 6 & 9 & 20 & y\\
2 & 2 & 3 & 3 & 10 & 20 & y\\
1 & 2 & 3 & 4 & 10 & 20 & y\\
1 & 1 & 4 & 4 & 10 & 20 & y\\
1 & 2 & 2 & 5 & 10 & 20 & y\\
1 & 1 & 3 & 5 & 10 & 20 & y\\
1 & 1 & 2 & 6 & 10 & 20 & y\\
 &  &  &  &  &  & \\
\hline\hline \end{tabular}}\hfill\\[3mm]
\hfil Table III: Weights for $l=5$ and $d\le 20$

\noindent
\hbox{\begin{tabular}{||rrrrr|r|r||} \hline\hline
$n_1$  & $n_2$  & $n_3$  & $n_4$  & $n_5$  & d & $\IP^4$   \\ \hline
\del
35 & 38 & 400 & 946 & 1419 & 2838 & n \\
36 & 37 & 401 & 948 & 1422 & 2844 & n \\
33 & 35 & 408 & 952 & 1428 & 2856 & n \\
31 & 37 & 408 & 952 & 1428 & 2856 & n \\
29 & 39 & 408 & 952 & 1428 & 2856 & n \\
27 & 41 & 408 & 952 & 1428 & 2856 & n \\
33 & 41 & 403 & 954 & 1431 & 2862 & n \\
\enddel
35 & 39 & 405 & 958 & 1437 & 2874 & n \\
34 & 35 & 414 & 966 & 1449 & 2898 & n \\
32 & 37 & 414 & 966 & 1449 & 2898 & n \\
31 & 38 & 414 & 966 & 1449 & 2898 & n \\
29 & 40 & 414 & 966 & 1449 & 2898 & n \\
28 & 41 & 414 & 966 & 1449 & 2898 & y \\
34 & 41 & 409 & 968 & 1452 & 2904 & n \\
33 & 37 & 420 & 980 & 1470 & 2940 & n \\
31 & 39 & 420 & 980 & 1470 & 2940 & n \\
29 & 41 & 420 & 980 & 1470 & 2940 & n \\
35 & 41 & 415 & 982 & 1473 & 2946 & n \\
35 & 36 & 426 & 994 & 1491 & 2982 & n \\
34 & 37 & 426 & 994 & 1491 & 2982 & n \\
33 & 38 & 426 & 994 & 1491 & 2982 & n \\
32 & 39 & 426 & 994 & 1491 & 2982 & n \\
31 & 40 & 426 & 994 & 1491 & 2982 & n \\
30 & 41 & 426 & 994 & 1491 & 2982 & n \\
29 & 42 & 426 & 994 & 1491 & 2982 & n \\
36 & 41 & 421 & 996 & 1494 & 2988 & y \\
35 & 37 & 432 & 1008 & 1512 & 3024 & n \\
31 & 41 & 432 & 1008 & 1512 & 3024 & n \\
36 & 37 & 438 & 1022 & 1533 & 3066 & n \\
35 & 38 & 438 & 1022 & 1533 & 3066 & n \\
34 & 39 & 438 & 1022 & 1533 & 3066 & n \\
33 & 40 & 438 & 1022 & 1533 & 3066 & n \\
32 & 41 & 438 & 1022 & 1533 & 3066 & n \\
31 & 42 & 438 & 1022 & 1533 & 3066 & n \\
35 & 39 & 444 & 1036 & 1554 & 3108 & n \\
33 & 41 & 444 & 1036 & 1554 & 3108 & n \\
37 & 38 & 450 & 1050 & 1575 & 3150 & n \\
34 & 41 & 450 & 1050 & 1575 & 3150 & n \\
37 & 39 & 456 & 1064 & 1596 & 3192 & n \\
35 & 41 & 456 & 1064 & 1596 & 3192 & n \\
38 & 39 & 462 & 1078 & 1617 & 3234 & n \\
37 & 40 & 462 & 1078 & 1617 & 3234 & n \\
36 & 41 & 462 & 1078 & 1617 & 3234 & y \\
37 & 41 & 468 & 1092 & 1638 & 3276 & n \\
39 & 40 & 474 & 1106 & 1659 & 3318 & n \\
38 & 41 & 474 & 1106 & 1659 & 3318 & n \\
37 & 42 & 474 & 1106 & 1659 & 3318 & n \\
39 & 41 & 480 & 1120 & 1680 & 3360 & n \\
40 & 41 & 486 & 1134 & 1701 & 3402 & n \\
41 & 42 & 498 & 1162 & 1743 & 3486 & y \\
\hline\hline \end{tabular}}\hfill\
\hbox{\begin{tabular}{||rrrrr|r|r||} \hline\hline
$n_1$  & $n_2$  & $n_3$  & $n_4$  & $n_5$  & d & $\IP^4$   \\ \hline
\del
20 & 33 & 199 & 471 & 690 & 1413 & n \\
19 & 35 & 200 & 473 & 692 & 1419 & n \\
13 & 42 & 204 & 476 & 693 & 1428 & n \\
15 & 38 & 204 & 476 & 695 & 1428 & n \\
21 & 32 & 201 & 476 & 698 & 1428 & n \\
19 & 30 & 204 & 476 & 699 & 1428 & n \\
21 & 26 & 204 & 476 & 701 & 1428 & n \\
\enddel
19 & 36 & 203 & 480 & 702 & 1440 & n \\
14 & 41 & 207 & 483 & 704 & 1449 & y \\
16 & 37 & 207 & 483 & 706 & 1449 & n \\
17 & 35 & 207 & 483 & 707 & 1449 & n \\
19 & 31 & 207 & 483 & 709 & 1449 & n \\
20 & 29 & 207 & 483 & 710 & 1449 & n \\
17 & 36 & 210 & 490 & 717 & 1470 & n \\
21 & 34 & 207 & 490 & 718 & 1470 & n \\
19 & 32 & 210 & 490 & 719 & 1470 & n \\
15 & 41 & 213 & 497 & 725 & 1491 & n \\
16 & 39 & 213 & 497 & 726 & 1491 & n \\
17 & 37 & 213 & 497 & 727 & 1491 & n \\
18 & 35 & 213 & 497 & 728 & 1491 & n \\
19 & 33 & 213 & 497 & 729 & 1491 & n \\
20 & 31 & 213 & 497 & 730 & 1491 & n \\
21 & 29 & 213 & 497 & 731 & 1491 & n \\
17 & 38 & 216 & 504 & 737 & 1512 & n \\
19 & 34 & 216 & 504 & 739 & 1512 & n \\
16 & 41 & 219 & 511 & 746 & 1533 & n \\
17 & 39 & 219 & 511 & 747 & 1533 & n \\
18 & 37 & 219 & 511 & 748 & 1533 & n \\
19 & 35 & 219 & 511 & 749 & 1533 & n \\
20 & 33 & 219 & 511 & 750 & 1533 & n \\
21 & 31 & 219 & 511 & 751 & 1533 & n \\
17 & 40 & 222 & 518 & 757 & 1554 & n \\
19 & 36 & 222 & 518 & 759 & 1554 & n \\
21 & 32 & 222 & 518 & 761 & 1554 & n \\
17 & 41 & 225 & 525 & 767 & 1575 & n \\
19 & 37 & 225 & 525 & 769 & 1575 & n \\
17 & 42 & 228 & 532 & 777 & 1596 & n \\
21 & 34 & 228 & 532 & 781 & 1596 & n \\
18 & 41 & 231 & 539 & 788 & 1617 & y \\
19 & 39 & 231 & 539 & 789 & 1617 & n \\
20 & 37 & 231 & 539 & 790 & 1617 & n \\
19 & 40 & 234 & 546 & 799 & 1638 & n \\
19 & 41 & 237 & 553 & 809 & 1659 & n \\
20 & 39 & 237 & 553 & 810 & 1659 & n \\
21 & 37 & 237 & 553 & 811 & 1659 & n \\
19 & 42 & 240 & 560 & 819 & 1680 & n \\
21 & 38 & 240 & 560 & 821 & 1680 & n \\
20 & 41 & 243 & 567 & 830 & 1701 & n \\
21 & 40 & 246 & 574 & 841 & 1722 & n \\
21 & 41 & 249 & 581 & 851 & 1743 & y \\
\hline\hline \end{tabular}}\hfill\\[3mm]
\hfil Table IV: Weights for $l=5$ and large $d$ (for $q_5=1/2$ and $q_5<1/2$)


\newpage

\def\LLab#1{\BP(0,0)\unitlength=1mm\put(-15,0){\makebox(0,0)[br]{\small#1}}\EP}
\def\ifundefined#1{\expandafter\ifx\csname#1\endcsname\relax}


							\ifundefined{draftmode}

\end{document}